\documentclass[11pt]{article}
\usepackage{latexsym}
\newcommand{\eproof}{\mbox{\ }\hfill $\Box$ \par \vskip 10pt}

\newtheorem{Theorem}{Theorem}[section]
\newtheorem{lemma}[Theorem]{Lemma}
\newtheorem{prop}[Theorem]{Proposition}

\begin{document}

\title{Dispersive estimates for the Schr\"odinger equation in dimensions four and five}

\author{{\sc Fernando Cardoso, Claudio Cuevas and Georgi Vodev}}

\date{}

\maketitle

\abstract{We prove optimal (that is, without loss of derivatives) dispersive estimates for the Schr\"odinger group
$e^{it(-\Delta+V)}$ for a class of real-valued potentials $V\in C^k({\bf R}^n)$, $V(x)=O(\langle x\rangle^{-\delta})$, where $n=4,5$, $k>(n-3)/2$, $\delta>3$ if $n=4$ and $\delta>5$ if $n=5$.}\\

\setcounter{section}{0}
\section{Introduction and statement of results}

The problem we address in the present paper is that one of finding the smallest possible regularity of a potential so that the corresponding Schr\"odinger group satisfies $L^1\to L^\infty$ dispersive estimates without loss of derivatives. While in dimensions one, two and three no regularity of the potential is required in order that such estimates hold true (see \cite{kn:G}, \cite{kn:GS}, \cite{kn:M}, \cite{kn:S}, \cite{kn:RS},  \cite{kn:V1}), the problem turns out to be quite hard in higher dimensions and far from being solved. 
Denote by $G$ the self-adjoint realization of the
operator $-\Delta+V$ on $L^2({\bf R}^n)$, $n\ge 4$, where $V\in L^\infty({\bf R}^n)$ is a real-valued potential satisfying 
$$|V(x)|\le C\langle x\rangle^{-\delta},\quad\forall x\in {\bf
R}^n,\eqno{(1.1)}$$ with constants $C>0$, $\delta>(n+2)/2$. It is well known that $G$ has no strictly positive eigenvalues and resonances 
similarly to the self-adjoint realization, $G_0$, of the free Laplacian $-\Delta$ on $L^2({\bf R}^n)$. The operator $G$, however, may have,
in contrast to the operator $G_0$, a finite number of non-positive eigenvalues and zero may be a resonance. It is shown in \cite{kn:V2} that under (1.1) only, the operator $e^{itG}\chi_a(G)$, $\forall a>0$, satisfies $L^1\to L^\infty$ dispersive estimates with a loss of $(n-3)/2$ derivatives, where $\chi_a\in C^\infty((-\infty,+\infty))$,
$\chi_a(\lambda)=0$ for $\lambda\le a$, $\chi_a(\lambda)=1$ for $\lambda\ge 2a$. On the other hand, the counterexample of \cite{kn:GoV}
shows the existence of potentials $V\in C^k_0({\bf R}^n)$, $\forall k<(n-3)/2$, for which $e^{itG}\chi_a(G)$ does not satisfy optimal
(that is, without loss of derivatives) dispersive estimates. Roughly speaking, the minimum regularity of the potential in order that optimal dispersive estimates could hold true is $V\in C^{(n-3)/2}({\bf R}^n)$. The problem of proving such dispersive estimates for potentials with this minimal regularity (even if they are supposed of compact support and small), however, is difficult and, to our best knowledege, is still open. Note that optimal dispersive estimates for $e^{itG}\chi_a(G)$ are proved in \cite{kn:JSS} for potentials satisfying (1.1) with $\delta>n$ as well as the condition
$$\widehat V\in L^1.\eqno{(1.2)}$$ 
This result has been recently extended in \cite{kn:MV} to potentials satisfying (1.1) with $\delta>n-1$ and (1.2). Clearly, (1.2) is fulfilled for potentials belonging to the Sobolev space $H^{n/2+\epsilon}({\bf R}^n)$, $\epsilon>0$. 
Roughly speaking, (1.2) requires $3/2+\epsilon$ more derivatives on the potential than what the counterexample of \cite{kn:GoV} would suggest in order to have optimal dispersive estimates. In the present paper we show that this gap can be reduced significantly when $n=4,5$. The problem, however, remains open when $n\ge 6$. To be more precise, given any $k,\delta\ge 0$, we introduce the space ${\cal C}^k_\delta({\bf R}^n)$ of
all functions $V\in C^k({\bf R}^n)$ satisfying
$$\left\|V\right\|_{{\cal C}_\delta^k}:=\sup_{x\in{\bf R}^n}\sum_{0\le|\alpha|\le k_0}\langle x\rangle^\delta\left|\partial_x^\alpha V(x)\right|$$ $$+\nu\sup_{x\in{\bf R}^n}\sum_{|\beta|=k_0}\langle x\rangle^\delta\sup_{x'\in{\bf R}^n:|x-x'|\le 1}\frac{\left|\partial_x^\beta V(x)-\partial_x^\beta V(x')\right|}{|x-x'|^\nu}<+\infty,$$
where $k_0\ge 0$ is an integer and $\nu=k-k_0$ satisfies $0\le\nu<1$. Introduce also 
the space ${\cal V}_\delta^k({\bf R}^n)$ of all functions $V\in C^k({\bf R}^n)$ satisfying
$$\left\|V\right\|_{{\cal V}_\delta^k}:=\sup_{x\in{\bf R}^n}\sum_{0\le|\alpha|\le k_0}\langle x\rangle^{\delta+|\alpha|}\left|\partial_x^\alpha V(x)\right|$$ $$+\nu\sup_{x\in{\bf R}^n}\sum_{|\beta|=k_0}\langle x\rangle^{\delta+k_0+1}\sup_{x'\in{\bf R}^n:|x-x'|\le 1}\frac{\left|\partial_x^\beta V(x)-\partial_x^\beta V(x')\right|}{|x-x'|^\nu}<+\infty,$$
where $k_0$ and $\nu$ are as above. Our main result is the following

\begin{Theorem} Let $n=4$ or $n=5$ and let $V\in{\cal C}_\delta^k({\bf R}^n)$ with $k>(n-3)/2$, $\delta>3$ if $n=4$ and $\delta>5$ if $n=5$. Then, for every 
$a>0$ there exists a constant $C=C(a)>0$ so that we have the dispersive estimate
$$\left\|e^{itG}\chi_a(G)\right\|_{L^1\to L^\infty}\le C|t|^{-n/2},\quad t\neq 0.\eqno{(1.3)}$$ 
Moreover, if in addition zero is neither an eigenvalue nor a resonance, then we have
$$\left\|e^{itG}P_{ac}\right\|_{L^1\to L^\infty}\le C|t|^{-n/2},\quad t\neq 0,\eqno{(1.4)}$$ 
where $P_{ac}$ denotes the spectral projection onto the absolutely continuous spectrum of $G$.
\end{Theorem}

Note that (1.4) follows from (1.3) and the low-frequency dispersive estimates proved in \cite{kn:MV} for potentials satisfying (1.1).

We expect that (1.3) holds true in any dimension $n\ge 4$ for potentials $V\in{\cal C}_\delta^k({\bf R}^n)$ with $k>(n-3)/2$, $\delta>n-1$. In fact, as suggested by the important Proposition 2.5 below, the condition (1.2) could probably be replaced by the following one
$$V\in{\cal V}_\delta^k({\bf R}^n) \quad\mbox{with}\quad k>(n-3)/2,\,\delta>2.\eqno{(1.5)}$$
In other words, there should be optimal dispersive estimates for potentials satisfying (1.5) as well as (1.1) with $\delta>n-1$ (or probably under (1.5) only). 

To prove (1.3) we follow the same strategy as in Appendix B of \cite{kn:MV} where optimal dispersive estimates have been proved for potentials satisfying (1.2) as well as (1.1) with $\delta>n-1$. The analysis, however, becomes harder without (1.2). Nevertheless, many arguments from \cite{kn:MV} still work in our case. The advantage of considering the case $n=4,5$ is that we can reduce the problem to estimating the $L^1\to L^\infty$ norm of operators (denoted by $T_1$ and $T_2$ below) with explicit kernels (see Proposition 2.2 below). These operators are obtained by iterating twice the {\it semi-classical} Duhamel formula (see (2.9)). In higher dimensions one must iterate this formula a finite number of times (depending on the dimension) and hence one must involve more operators most of which being no longer with explicit kernels (see \cite{kn:CV} for more details). In other words, the bigger the dimension is, the more difficult the proof is. To make our approach work in higher dimensions we also need to improve in $h$ the estimates (2.16) and (3.1) below concerning the operator $T_2$ provided more regularity on the potential is required. This turns out to be hard despite the fact that $T_2$ has an explicit kernel. In contrast, the operator $T_1$ can be treated in all dimensions (see Propositions 2.5 and 2.6 below).\\

{\bf Acknowledgements.} The authors have been supported by the agreement Brazil-France in Mathematics - Proc. 69.0014/01-5. The first two authors have also been  partially supported by the CNPq-Brazil.

\section{Proof of Theorem 1.1} Without loss of generality we may suppose $t>0$. Given a function $\psi\in C_0^\infty((0,+\infty))$ and a
parameter $0<h\le 1$, following \cite{kn:V1}, \cite{kn:V2} (see also Appendix B of \cite{kn:MV}), we set
$$\Psi(t,h)=e^{itG}\psi(h^2G)-e^{itG_0}\psi(h^2G_0),$$
$$F(t)=i\int_0^te^{i(t-\tau)G_0}Ve^{i\tau G_0}d\tau.$$
By a simple argument originating from \cite{kn:V1} (see also \cite{kn:V2}) it is easy to see that (1.3) follows from the following
{\it semi-classical} dispersive estimates.

\begin{Theorem} Under the assumptions of Theorem 1.1, there exist constants $C,\beta>0$ so that for all $t>0$, $0<h\le 1$, we have the estimates
$$\left\|F(t)\right\|_{L^1\to L^\infty}\le Ct^{-n/2},\eqno{(2.1)}$$
$$\left\|\Psi(t,h)-F(t)\psi(h^2G_0)\right\|_{L^1\to L^\infty}\le Ch^\beta t^{-n/2}.\eqno{(2.2)}$$
\end{Theorem}

{\it Proof.} A lot of our analysis works out in all dimensions $n\ge 4$, so in what follows we will establish our key estimates in the greatest possible generality. This might be usefull in view of possible further applications. Let $\psi_1\in C_0^\infty((0,+\infty))$ be such that $\psi_1=1$ on supp$\,\psi$. Define the operators $T_j(t,h)$, $j=1,2$, as follows
$$T_1(t,h)=i\int_0^te^{i(t-\tau)G_0}\psi_1(h^2G_0)Ve^{i\tau G_0}\psi(h^2G_0)d\tau,$$
$$T_2(t,h)=i\int_0^te^{i(t-\tau)G_0}\psi_1(h^2G_0)VT_1(\tau,h)d\tau.$$
Using some estimates from \cite{kn:V2} we will first prove the following 

\begin{prop} Let $V$ satisfy (1.1). Then, for all $t>0$, $0<h\le 1$, $0<\epsilon\ll 1$, we have the estimates
$$\left\|T_2(t,h)\right\|_{L^1\to L^\infty}\le Ch^{-(n-4)/2} t^{-n/2},\eqno{(2.3)}$$
$$\left\|\Psi(t,h)-\sum_{j=1}^2T_j(t,h)\right\|_{L^1\to L^\infty}\le Ch^{-(n-6)/2-\epsilon} t^{-n/2}.\eqno{(2.4)}$$
\end{prop} 

{\it Proof.} To prove (2.3) we will use the following estimates proved in \cite{kn:V2} (see Proposition 2.1):
$$\left\|\langle x\rangle^{-s}e^{it G_0}\psi(h^2G_0)\langle x\rangle^{-s}\right\|_{L^2\to L^2}\le C\langle t/h\rangle^{-s},\quad\forall t,\,
0<h\le 1,\,\,s\ge 0,\eqno{(2.5)}$$
$$\left\|e^{it G_0}\psi(h^2G_0)\langle x\rangle^{-1/2-s-\epsilon}\right\|_{L^2\to L^\infty}\le Ch^{s-(n-1)/2}|t|^{-s-1/2},\eqno{(2.6)}$$
$\forall t\neq 0,\,0<h\le 1,\,\,0\le s\le (n-1)/2$. For $1/2-\epsilon/4\le s\le (n-1)/2$, $0<\epsilon\ll 1$, we get
$$t^{s+1/2}\left\|\langle x\rangle^{-1/2-s-\epsilon}T_1(t,h)\right\|_{L^1\to L^2}$$ $$\le C\int_0^{t/2}(t-\tau)^{s+1/2}
\left\|\langle x\rangle^{-1/2-s-\epsilon}e^{i(t-\tau) G_0}\psi_1(h^2G_0)\langle x\rangle^{-1/2-s-\epsilon}\right\|_{L^2\to L^2}$$ $$\times 
\left\|\langle x\rangle^{-1-\epsilon}e^{i\tau G_0}\psi(h^2G_0)\right\|_{L^1\to L^2}d\tau$$
 $$+C\int_{t/2}^t\tau^{s+1/2}
\left\|\langle x\rangle^{-1/2-s-\epsilon}e^{i(t-\tau) G_0}\psi_1(h^2G_0)\langle x\rangle^{-1-\epsilon}\right\|_{L^2\to L^2}$$ $$\times 
\left\|\langle x\rangle^{-1/2-s-\epsilon}e^{i\tau G_0}\psi(h^2G_0)\right\|_{L^1\to L^2}d\tau$$
$$\le Ch^{s+1/2}\int_0^\infty\left\|\langle x\rangle^{-1-\epsilon}e^{i\tau G_0}\psi(h^2G_0)\right\|_{L^1\to L^2}d\tau$$
$$+Ch^{s-(n-1)/2}\int_0^\infty\left\|\langle x\rangle^{-1-\epsilon/2}e^{i\tau G_0}\psi_1(h^2G_0)\langle x\rangle^{-1-\epsilon/2}\right\|_{L^2\to L^2}d\tau$$
 $$\le Ch^{s-(n-3)/2+\beta'}\int_h^\infty\tau^{-1-\beta'}d\tau+Ch^{s-(n-3)/2-\beta'}\int_0^h\tau^{-1+\beta'}d\tau$$
$$+Ch^{s-(n-1)/2}\int_0^\infty\langle\tau/h\rangle^{-1-\epsilon/2}d\tau\le Ch^{s-(n-3)/2},\eqno{(2.7)}$$
where $0<\beta'\ll 1$ is independent of $h$. In the same way, using (2.6) and (2.7), we obtain
$$t^{n/2}\left\|T_2(t,h)\right\|_{L^1\to L^\infty}$$ $$\le C\int_0^{t/2}(t-\tau)^{n/2}
\left\|e^{i(t-\tau) G_0}\psi_1(h^2G_0)\langle x\rangle^{-n/2-\epsilon}\right\|_{L^2\to L^\infty} 
\left\|\langle x\rangle^{-1-\epsilon}T_1(\tau,h)\right\|_{L^1\to L^2}d\tau$$
 $$+C\int_{t/2}^t\tau^{n/2}
\left\|e^{i(t-\tau) G_0}\psi_1(h^2G_0)\langle x\rangle^{-1-\epsilon}\right\|_{L^2\to L^\infty}
\left\|\langle x\rangle^{-n/2-\epsilon}T_1(\tau,h)\right\|_{L^1\to L^2}d\tau$$
 $$\le C\int_0^\infty
\left\|\langle x\rangle^{-1-\epsilon}T_1(\tau,h)\right\|_{L^1\to L^2}d\tau+Ch\int_0^\infty
\left\|e^{i\tau G_0}\psi_1(h^2G_0)\langle x\rangle^{-1-\epsilon}\right\|_{L^2\to L^\infty}d\tau$$ $$\le Ch^{-(n-4)/2}.\eqno{(2.8)}$$
To prove (2.4) we will make use of the following identity which can be derived easily from Duhamel's formula (see \cite{kn:V1}, \cite{kn:V2}):
$$\Psi(t,h)=\sum_{j=1}^3\Psi_j(t,h),\eqno{(2.9)}$$
where
$$\Psi_1(t,h)=\left(\psi_1(h^2G)-\psi_1(h^2G_0)\right)\Psi(t,h)+\psi_1(h^2G_0)e^{it G_0}\left(\psi(h^2G)-\psi(h^2G_0)\right)$$
$$+\left(\psi_1(h^2G)-\psi_1(h^2G_0)\right)e^{it G_0}\psi(h^2G_0),$$
$$\Psi_2(t,h)=T_1(t,h),$$
$$\Psi_3(t,h)=i\int_0^te^{i(t-\tau)G_0}\psi_1(h^2G_0)V\Psi(\tau,h)d\tau.$$
Iterating this identity once again, we obtain
$$\Psi_3(t,h)=\sum_{j=1}^3\Psi_3^{(j)}(t,h),\eqno{(2.10)}$$
where
$$\Psi_3^{(1)}(t,h)=T_2(t,h),$$
$$\Psi_3^{(2)}(t,h)=i\int_0^te^{i(t-\tau)G_0}\psi_1(h^2G_0)V\Psi_1(\tau,h)d\tau,$$
$$\Psi_3^{(3)}(t,h)=i\int_0^t\widetilde T_1(t-\tau,h)V\Psi(\tau,h)d\tau,$$
where $\widetilde T_1$ is defined by replacing in the definition of $T_1$ the operator $\psi(h^2G_0)$ by $\psi_1(h^2G_0)$.
Clearly, $\widetilde T_1$ satisfies (2.7). 
We will use now the following estimates proved in \cite{kn:V2}.

\begin{prop} Let $V$ satisfy (1.1). Then, for all $t\neq 0$, $0<h\le 1$, $0<\epsilon\ll 1$, $1/2-\epsilon/4\le s\le (n-1)/2$, we have the estimates
$$\left\|\langle x\rangle^{-1/2-s-\epsilon}\Psi(t,h)\right\|_{L^1\to L^2}\le Ch^{s-(n-3)/2-\epsilon} |t|^{-s-1/2},\eqno{(2.11)}$$
$$\left\|\Psi(t,h)\right\|_{L^1\to L^\infty}\le Ch^{-(n-3)/2}|t|^{-n/2},\eqno{(2.12)}$$
$$\left\|\langle x\rangle^{-1/2-s-\epsilon}\Psi_1(t,h)\right\|_{L^1\to L^2}\le Ch^{s-(n-7)/2-\epsilon}|t|^{-s-1/2},\eqno{(2.13)}$$
$$\left\|\Psi_1(t,h)\right\|_{L^1\to L^\infty}\le Ch^{-(n-7)/2}|t|^{-n/2}.\eqno{(2.14)}$$
\end{prop}

Combining (2.11) and (2.13) together with (2.6) and (2.7), in the same way as above, it is easy to get the estimates
$$\left\|\Psi_3^{(j)}(t,h)\right\|_{L^1\to L^\infty}\le Ch^{-(n-6)/2-\epsilon} t^{-n/2},\quad j=2,3.\eqno{(2.15)}$$
Now (2.4) follows from (2.9), (2.10), (2.14) and (2.15).
\eproof

Thus the problem of proving (2.2) is reduced to studying the operators $T_1$ and $T_2$. To make our proof work when $n=4$ we need to improve the estimate (2.3) with an extra factor $O(h^\epsilon)$. This can be done if one requires a little regularity on the potential. More precisely, we have the following

\begin{prop} Let $V\in{\cal C}_\delta^\varepsilon({\bf R}^n)$ with $0<\varepsilon\ll 1$, $\delta>(n+2)/2$. Then, there exist constants
$C,\varepsilon_0>0$ so that for all $t>0$, $0<h\le 1$, we have the estimate
$$\left\|T_2(t,h)\right\|_{L^1\to L^\infty}\le Ch^{-(n-4)/2+\varepsilon_0} t^{-n/2}.\eqno{(2.16)}$$
\end{prop}

{\it Proof.} Let $\rho\in C_0^\infty({\bf R}^n)$ be a real-valued function such that $\int\rho(x)dx=1$, and set $\rho_\theta(x)=\theta^{-n}\rho(x/\theta)$, $0<\theta\le 1$. It is easy to see that $V\in{\cal C}_\delta^\varepsilon({\bf R}^n)$ implies that the function $V_\theta=V*\rho_\theta$ satisfies the bounds
$$\left|V_\theta(x)\right|\le C\langle x\rangle^{-\delta},\quad\forall x\in {\bf R}^n,\eqno{(2.17)}$$
$$\left|V(x)-V_\theta(x)\right|\le C\theta^\varepsilon\langle x\rangle^{-\delta},\quad\forall x\in {\bf R}^n,\eqno{(2.18)}$$
with a constant $C>0$ independent of $\theta$. We also have
$$\left\|\widehat V_\theta\right\|_{L^1}=\left\|\widehat\rho_\theta\widehat V\right\|_{L^1}\le \left\|\widehat\rho_\theta\right\|_{L^2}\left\|\widehat V\right\|_{L^2}=Const\left\|\rho_\theta\right\|_{L^2}\left\|V\right\|_{L^2}\le C\theta^{-n/2}.$$
Therefore, we get (see \cite{kn:JSS})
$$\left\|e^{-itG_0} V_\theta e^{itG_0}\right\|_{L^1\to L^1}\le \left\|\widehat V_\theta\right\|_{L^1}\le C\theta^{-n/2},\quad\forall t.\eqno{(2.19)}$$
Define now the operator $T_{2,\theta}(t,h)$ by replacing in the definition of 
$T_2(t,h)$ the potential $V$ by $V_\theta$. In the same way as in (2.8), using (2.17) and (2.18), we get
$$\left\|T_2(t,h)-T_{2,\theta}(t,h)\right\|_{L^1\to L^\infty}\le C\theta^\varepsilon h^{-(n-4)/2} t^{-n/2}.\eqno{(2.20)}$$
Let $0<\gamma<1$ be a parameter to be fixed later on depending on $h$. For $t\ge 4\gamma$, we decompose the operator $T_{2,\theta}$ as $T_{2,\theta}^{(1)}+T_{2,\theta}^{(2)}$, where
$$T_{2,\theta}^{(j)}(t,h)=i^2\int\int_{I_j}e^{i(t-\tau)G_0}\psi_1(h^2G_0)V_\theta e^{i(\tau-s)G_0}\psi_1(h^2G_0)V_\theta e^{isG_0}\psi(h^2G_0)dsd\tau,$$
$$I_1=[\gamma,t-\gamma]\times[0,\tau]\cup[t-\gamma,t]\times[\gamma,\tau-\gamma],$$
$$I_2=[0,\gamma]\times[0,\tau]\cup[t-\gamma,t]\times[0,\gamma]\cup[t-\gamma,t]\times[\tau-\gamma,\tau].$$
Proceeding as in (2.8) and using (2.17) one can easily get the estimate
$$\left\|T_{2,\theta}^{(1)}(t,h)\right\|_{L^1\to L^\infty}\le C\gamma^{-\epsilon} h^{-(n-4)/2+\epsilon} t^{-n/2}\eqno{(2.21)}$$
for every $0<\epsilon\ll 1$. On the other hand, using (2.19) together with the fact that the operator $\psi(h^2G_0)$ is uniformly bounded on $L^1$, we obtain
$$\left\|T_{2,\theta}^{(2)}(t,h)\right\|_{L^1\to L^\infty}\le C\gamma^2\theta^{-n}t^{-n/2},\eqno{(2.22)}$$
for $t\ge 4\gamma$. Clearly, for $0<t\le 4\gamma$ the estimate (2.22) holds with $T_{2,\theta}^{(2)}$ replaced by $T_{2,\theta}$. 
Combining (2.20)-(2.22) we conclude
$$\left\|T_2(t,h)\right\|_{L^1\to L^\infty}\le C\left(\theta^\varepsilon+\gamma^{-\epsilon} h^{\epsilon}+\gamma^2\theta^{-n}\right) h^{-(n-4)/2} t^{-n/2}.\eqno{(2.23)}$$
Now, taking $\gamma=h^{1/2}$, $\theta=h^{1/(n+\varepsilon)}$, we deduce (2.16) from (2.23).
\eproof

To deal with the operator $T_1$ we need the following proposition the proof of which will be given in Section 4.

\begin{prop} Let $V\in{\cal V}_\delta^{(n-3)/2+\varepsilon}({\bf R}^n)$ with $0<\varepsilon\ll 1$, $\delta>2$. Then, there exist constants
$C>0$, $0<\varepsilon'\ll 1$, so that for every $0<\gamma\le 1$ we have the estimates
$$\left\|\int_0^\gamma e^{i(t-\tau)G_0}Ve^{i\tau G_0}d\tau\right\|_{L^1\to L^\infty}\le C\gamma^{\varepsilon'} t^{-n/2},\quad t\ge 2\gamma,\eqno{(2.24)}$$
$$\left\|\int_0^te^{i(t-\tau)G_0}Ve^{i\tau G_0}d\tau\right\|_{L^1\to L^\infty}\le Ct^{-n/2+\varepsilon'},\quad 0<t\le 2.\eqno{(2.25)}$$
\end{prop}

Clearly, for $0<t\le 2$ the estimate (2.1) follows from (2.25). On the other hand, it is shown in Appendix B of \cite{kn:MV} that for $t\ge 2$ we have the estimate
$$\left\|\int_1^{t-1}e^{i(t-\tau)G_0}Ve^{i\tau G_0}d\tau\right\|_{L^1\to L^\infty}\le Ct^{-n/2}\eqno{(2.26)}$$
for potentials satisfying (1.1) with $\delta>n-1$. So, for $t\ge 2$ the estimate (2.1) follows from (2.24) and (2.26). It is also clear that in the particular case of $n=4$ the estimate (2.2) follows from Propositions 2.2, 2.4 and the following

\begin{prop} Let $V\in{\cal V}_\delta^{(n-3)/2+\varepsilon}({\bf R}^n)$ with $0<\varepsilon\ll 1$, $\delta>2$, and suppose in addition that $V$ satisfies (1.1) with $\delta>n-1$. Then, there exist constants
$C>0$, $0<\varepsilon_1\ll 1$, so that we have the estimate
$$\left\|\int_0^t(1-\psi_1)(h^2G_0)e^{i(t-\tau)G_0}Ve^{i\tau G_0}\psi(h^2G_0)d\tau\right\|_{L^1\to L^\infty}\le Ch^{\varepsilon_1}t^{-n/2}.\eqno{(2.27)}$$
\end{prop}

{\it Proof.} Let $0<\gamma<1$ be a parameter to be fixed later on depending on $h$. By (2.24),
$$\left\|\left(\int_0^\gamma+\int_{t-\gamma}^t\right)(1-\psi_1)(h^2G_0)e^{i(t-\tau)G_0}Ve^{i\tau G_0}\psi(h^2G_0)d\tau\right\|_{L^1\to L^\infty}\le C\gamma^{\varepsilon'}t^{-n/2}.\eqno{(2.28)}$$
On the other hand, it is proved in \cite{kn:MV} (see Proposition B.5) that for potentials satisfying (1.1) with $\delta>n-1$ we have the estimate
$$\left\|\int_\gamma^{t-\gamma}(1-\psi_1)(h^2G_0)e^{i(t-\tau)G_0}Ve^{i\tau G_0}\psi(h^2G_0)d\tau\right\|_{L^1\to L^\infty}\le Ch^\epsilon\gamma^{-(n-3)/2-\epsilon}t^{-n/2},\eqno{(2.29)}$$
for every $0<\epsilon\ll 1$. By (2.28) and (2.29),
$$\left\|\int_0^t(1-\psi_1)(h^2G_0)e^{i(t-\tau)G_0}Ve^{i\tau G_0}\psi(h^2G_0)d\tau\right\|_{L^1\to L^\infty}\le C\left(\gamma^{\varepsilon'}+h^\epsilon\gamma^{-(n-3)/2-\epsilon}\right)t^{-n/2}.\eqno{(2.30)}$$
Now, taking $\gamma=h^{\beta_1}$ with a suitably chosen constant $\beta_1>0$, we get (2.27).
\eproof

To prove (2.2) in the case $n=5$ we need to improve in $h$ the estimate (2.16) provided more regularity of the potential is required. To do so, we have to refine the estimates in the proof of Proposition 2.4 and especially (2.22). This will be carried out in the next section.

\section{Study of the operator $T_2$ in the case $n=5$}

In this section we will prove the following

\begin{prop} Let $V\in{\cal C}_\delta^1({\bf R}^5)$ with $\delta>5$. Then, there exist constants
$C,\varepsilon_0>0$ so that for all $t>0$, $0<h\le 1$, we have the estimate
$$\left\|T_2(t,h)\right\|_{L^1\to L^\infty}\le Ch^{\varepsilon_0} t^{-5/2}.\eqno{(3.1)}$$
\end{prop}

{\it Proof.} We keep the same notations as in the previous section. It is easy to see that $V\in{\cal C}_\delta^1({\bf R}^5)$ implies that the function $V_\theta=V*\rho_\theta$ satisfies the bounds
$$\left|V_\theta(x)\right|\le C\langle x\rangle^{-\delta},\quad\forall x\in {\bf R}^5,\eqno{(3.2)}$$
$$\left|V(x)-V_\theta(x)\right|\le C\theta\langle x\rangle^{-\delta},\quad\forall x\in {\bf R}^5,\eqno{(3.3)}$$
$$\left|\partial_x^\alpha V_\theta(x)\right|\le C_\alpha\theta^{1-|\alpha|}\langle x\rangle^{-\delta},\quad\forall x\in {\bf R}^5,\quad |\alpha|\ge 1,\eqno{(3.4)}$$
with constants $C,C_\alpha>0$ independent of $\theta$. In view of (3.4) we have
$$\left\|\widehat V_\theta\right\|_{L^1}=\int_{|\xi|\le R}\left|\widehat V_\theta(\xi)\right|d\xi+\int_{|\xi|\ge R}\left|\widehat V_\theta(\xi)\right|d\xi$$
$$\le \left(\int|\xi|^4\left|\widehat V_\theta(\xi)\right|^2d\xi\right)^{1/2}\left(\int_{|\xi|\le R}|\xi|^{-4}d\xi\right)^{1/2}$$
 $$+\left(\int|\xi|^6\left|\widehat V_\theta(\xi)\right|^2d\xi\right)^{1/2}\left(\int_{|\xi|\ge R}|\xi|^{-6}d\xi\right)^{1/2}$$
 $$\le CR^{1/2}\sum_{|\alpha|=2}\left\|\partial_x^\alpha V_\theta(x)\right\|_{L^2({\bf R}^5)}+CR^{-1/2}\sum_{|\alpha|=3}\left\|\partial_x^\alpha V_\theta(x)\right\|_{L^2({\bf R}^5)}$$
 $$\le CR^{1/2}\theta^{-1}+CR^{-1/2}\theta^{-2}=O\left(\theta^{-3/2}\right),\eqno{(3.5)}$$
 if we choose $R=\theta^{-1}$.  We have the following

\begin{lemma} For all $t>0$, $0<h\le 1$, $0<\theta\le 1$, we have
$$\left\|T_2(t,h)-T_{2,\theta}(t,h)\right\|_{L^1\to L^\infty}\le C\theta h^{-1/2}t^{-5/2},\eqno{(3.6)}$$
with a constant $C>0$ independent of $t$, $h$ and $\theta$.
\end{lemma}

{\it Proof.} The estimate (3.6) is obtained in the same way as (2.20) using (3.3) instead of (2.18).
 \eproof
 
\begin{lemma} For all $t\ge 4\gamma$, $0<h\le 1$, $0<\theta\le 1$, $0<\gamma\le 1$, we have
$$\left\|T_{2,\theta}^{(1)}(t,h)\right\|_{L^1\to L^\infty}\le Ch\gamma^{-3/2}t^{-5/2},\eqno{(3.7)}$$
with a constant $C>0$ independent of $t$, $h$, $\theta$ and $\gamma$.
\end{lemma}

{\it Proof.} Using (2.6) (with $n=5$, $s=(n-1)/2$), (2.7) (with $V$ replaced by $V_\theta$, $n=5$, $s=(n-1)/2$) and (3.2) (with $\delta>5$) 
we get
$$t^{5/2}\left\|\int_\gamma^{t-\gamma}\int_0^\tau e^{i(t-\tau)G_0}\psi_1(h^2G_0)V_\theta e^{i(\tau-s)G_0}\psi_1(h^2G_0)V_\theta  e^{isG_0}\psi(h^2G_0)dsd\tau\right\|_{L^1\to L^\infty}$$ $$\le C\int_\gamma^{t/2}(t-\tau)^{5/2}
\left\|e^{i(t-\tau) G_0}\psi_1(h^2G_0)\langle x\rangle^{-5/2-\epsilon}\right\|_{L^2\to L^\infty}$$ $$\times  
\left\|\langle x\rangle^{-5/2-\epsilon}\int_0^\tau e^{i(\tau-s)G_0}\psi_1(h^2G_0)V_\theta e^{is G_0}\psi(h^2G_0)ds\right\|_{L^1\to L^2}d\tau$$
 $$+C\int_{t/2}^{t-\gamma}\tau^{5/2}
\left\|e^{i(t-\tau) G_0}\psi_1(h^2G_0)\langle x\rangle^{-5/2-\epsilon}\right\|_{L^2\to L^\infty}$$ $$\times 
\left\|\langle x\rangle^{-5/2-\epsilon}\int_0^\tau e^{i(\tau-s)G_0}\psi_1(h^2G_0)V_\theta e^{is G_0}\psi(h^2G_0)ds\right\|_{L^1\to L^2}d\tau$$
 $$\le Ch\int_\gamma^\infty\tau^{-5/2}d\tau\le Ch\gamma^{-3/2}.$$
To bound the norm of the integral over $[t-\gamma,t]\times[\gamma,\tau-\gamma]$ observe that it can be written in the form
$$\int_\gamma^{t-\gamma}M(t-s,h)V_\theta  e^{isG_0}\psi(h^2G_0)ds,$$
where the operator
$$M(t,h)=\int_0^\gamma e^{i\tau G_0}\psi_1(h^2G_0)V_\theta e^{i(t-\tau)G_0}\psi_1(h^2G_0)d\tau$$
satisfies the estimate
$$\left\|M(t,h)\langle x\rangle^{-5/2-\epsilon}\right\|_{L^2\to L^\infty}\le Cht^{-5/2}.$$
Therefore, it can be treated in the same way as above.
\eproof

\begin{prop} For all $t\ge 4\gamma$, $0<h\le 1$, $0<\theta\le 1$, $0<\gamma\le 1$, $0<\epsilon\ll 1$, we have
$$\left\|T_{2,\theta}^{(2)}(t,h)\right\|_{L^1\to L^\infty}\le C_\epsilon\gamma^{2-\epsilon}\theta^{-5/2}t^{-5/2},\eqno{(3.8)}$$
with a constant $C_\epsilon>0$ independent of $t$, $h$, $\theta$ and $\gamma$. Moreover, for $0<t\le 4\gamma$ the estimate (3.8) holds
with $T_{2,\theta}^{(2)}$ replaced by $T_{2,\theta}$. 
\end{prop}

{\it Proof.} In view of (3.5), we get 
$$\left\|e^{-itG_0} V_\theta e^{itG_0}\right\|_{L^1\to L^1}\le \left\|\widehat V_\theta\right\|_{L^1}\le C\theta^{-3/2},\quad\forall t.\eqno{(3.9)}$$
In what follows we will derive (3.8) from (3.9) and the following proposition the proof of which will be given in Section 5.

\begin{prop} Let $n$ be odd and let $V\in{\cal V}_\delta^{(n-1)/2}({\bf R}^n)$ with $\delta>1$. Then, for every $0<\gamma\le 1$, 
$0<\epsilon\ll 1$, we have the estimates
$$\left\|\int_0^\gamma e^{i(t-\tau)G_0}Ve^{i\tau G_0}d\tau\right\|_{L^1\to L^\infty}\le C_\epsilon\left\|V\right\|_{{\cal V}_\delta^{(n-1)/2}}\gamma^{1-\epsilon} t^{-n/2},\quad t\ge 2\gamma,\eqno{(3.10)}$$
$$\left\|\int_0^te^{i(t-\tau)G_0}Ve^{i\tau G_0}d\tau\right\|_{L^1\to L^\infty}\le C_\epsilon \left\|V\right\|_{{\cal V}_\delta^{(n-1)/2}}t^{-n/2+1-\epsilon},\quad 0<t\le 2,\eqno{(3.11)}$$
with a constant $C_\epsilon>0$ independent of $t$, $\gamma$ and $V$.
\end{prop}

We are going to use this proposition with $n=5$ and $V$ replaced by $V_\theta$. To this end, observe that (3.4) implies
$$\left\|V_\theta\right\|_{{\cal V}_\delta^2}\le C\theta^{-1}.\eqno{(3.12)}$$
Setting
$$P(\tau)=\int_0^\tau e^{i(t-\tau')G_0}V_\theta e^{i\tau' G_0}\psi_1(h^2G_0)d\tau',$$
$$Q(\tau)=\int_0^\tau e^{-isG_0}V_\theta e^{isG_0}ds,$$
we can write
$$\int_0^\gamma\int_0^\tau e^{i(t-\tau)G_0}V_\theta e^{i(\tau-s)G_0}\psi_1(h^2G_0)V_\theta  e^{isG_0}dsd\tau$$
$$=\int_0^\gamma P'(\tau)Q(\tau)d\tau=P(\gamma)Q(\gamma)-\int_0^\gamma P(\tau)Q'(\tau)d\tau.\eqno{(3.13)}$$
By Proposition 3.5 and (3.12) we have
$$\left\|P(\tau)\right\|_{L^1\to L^\infty}\le C_\epsilon\theta^{-1}\tau^{1-\epsilon}t^{-5/2}.\eqno{(3.14)}$$
On the other hand, (3.9) implies
$$\left\|Q'(\tau)\right\|_{L^1\to L^1}\le C\theta^{-3/2},\eqno{(3.15)}$$
$$\left\|Q(\tau)\right\|_{L^1\to L^1}\le C\tau\theta^{-3/2}.\eqno{(3.16)}$$
By (3.13)-(3.16) we conclude that the $L^1\to L^\infty$ norm of the LHS of (3.13) is upper bounded by
$$C_\epsilon\theta^{-5/2}\gamma^{2-\epsilon}t^{-5/2}.$$
It is easy to see that the $L^1\to L^\infty$ norm of the integrals over $[t-\gamma,t]\times[0,\gamma]$ and $[t-\gamma,t]\times[\tau-\gamma,\tau]$ can be bounded in the same way.
\eproof

By Lemmas 3.2 and 3.3 and Proposition 3.4, we conclude
$$\left\|T_2(t,h)\right\|_{L^1\to L^\infty}\le C\left(\theta h^{-1/2}+h\gamma^{-3/2}+\gamma^{2-\epsilon}\theta^{-5/2}\right) t^{-5/2}.\eqno{(3.17)}$$
Now, taking $\theta=h^{1/2+\epsilon}$, $\gamma=h^{2/3-\epsilon}$, we deduce (3.1) from (3.17).
\eproof

\section{Proof of Proposition 2.5}

It is easy to see that (2.25) follows from (2.24) applied with $\gamma=t/2$. To prove (2.24) observe that the kernel of the operator
$$\int_0^\gamma e^{i(t-\tau)G_0}Ve^{i\tau G_0}d\tau$$
is of the form
$$K(x,y,t,\gamma)=c_n\int_0^\gamma\int_{{\bf R}^n}e^{i\varphi}(t-\tau)^{-n/2}\tau^{-n/2}V(\xi)d\xi d\tau,$$
where $c_n$ is a constant and
$$\varphi(x,y,\xi,t,\tau)=\frac{|x-\xi|^2}{4(t-\tau)}+\frac{|y-\xi|^2}{4\tau}.$$
Let $\phi\in C_0^\infty({\bf R})$, $\phi(\lambda)=1$ for $|\lambda|\le 1/2$, $\phi(\lambda)=0$ for $|\lambda|\ge 1$. Clearly, the function $1-\phi$ can be written as
$$(1-\phi)(\lambda)=\sum_{q=1}^\infty\widetilde\phi(2^{-q}\lambda)$$
with a function $\widetilde\phi\in C_0^\infty({\bf R})$, $\widetilde\phi(\lambda)=0$ for $|\lambda|\le 1/2$ and $|\lambda|\ge 1$.
Therefore, we have
$$1=\sum_{q=0}^\infty\phi_q(\lambda),$$
where $\phi_0=\phi$, $\phi_q(\lambda)=\widetilde\phi(2^{-q}\lambda)$, $q\ge 1$. Write now the function $K$ as
$$K(x,y,t,\gamma)=c_n\sum_{p=1}^\infty\sum_{q=0}^\infty K_{p,q}(x,y,t,\gamma),\eqno{(4.1)}$$
where
$$K_{p,q}(x,y,t,\gamma)=\int_0^{\gamma}\int_{{\bf R}^n}e^{i\varphi}(t-\tau)^{-n/2}\tau^{-n/2}\phi_p(\gamma/\tau)\phi_q(|\xi|)V(\xi)d\xi d\tau.$$
Clearly, (2.24) would follow from (4.1) and the following

\begin{prop} Under the assumptions of Proposition 2.5, there exist constants $C>0$ and $0<\varepsilon'\ll 1$ so that we have the bound
$$\left|K_{p,q}(x,y,t,\gamma)\right|\le C2^{-\varepsilon'(p+q)}\gamma^{\varepsilon'} t^{-n/2},\quad t\ge 2\gamma.\eqno{(4.2)}$$
\end{prop}

{\it Proof.} Set $k_0=(n-4)/2$ if $n$ is even and 
$k_0=(n-3)/2$ if $n$ is odd, $k=(n-3)/2+\varepsilon$,  
$\nu=k-k_0$, $0<\nu<1$. Set also $V_q(\xi)=\phi_q(|\xi|)V(\xi)$, $V_{q,\theta}=\rho_\theta*V_q$, $0<\theta\le 1$, where the function $\rho_\theta$ is as in the previous section. Since $V\in {\cal V}^k_\delta({\bf R}^n)$, we have
$$\left|\partial_\xi^\alpha V_q(\xi)\right|\le C2^{-q(\delta+|\alpha|)},\quad 0\le|\alpha|\le k_0,\eqno{(4.3)}$$
$$\left|\partial_\xi^\alpha V_q(\xi)-\partial_\xi^\alpha V_q(\xi')\right|\le C2^{-q(\delta+|\alpha|+1)}|\xi-\xi'|^\nu,\quad|\xi-\xi'|\le 1,\quad |\alpha|= k_0.\eqno{(4.4)}$$
It is easy to see that these bounds imply
$$\left|\partial_\xi^\alpha V_{q,\theta}(\xi)\right|\le C2^{-q(\delta+|\alpha|)},\quad 0\le|\alpha|\le k_0,\eqno{(4.5)}$$
$$\left|\partial_\xi^\alpha V_{q,\theta}(\xi)\right|\le C\theta^{-1+\nu}2^{-q(\delta+|\alpha|)},\quad |\alpha|=k_0+1,\eqno{(4.6)}$$
$$\left|\partial_\xi^\alpha V_q(\xi)-\partial_\xi^\alpha V_{q,\theta}(\xi)\right|\le C\theta^\nu 2^{-q(\delta+|\alpha|)},\quad 0\le|\alpha|\le k_0.\eqno{(4.7)}$$
Decompose now the function $K_{p,q}$ as $K_{p,q}^{(1)}+K_{p,q}^{(2)}$, where $K_{p,q}^{(1)}$ and $K_{p,q}^{(2)}$ are defined by replacing $V_q$ in the definition of $K_{p,q}$ by $V_q-V_{q,\theta}$ and $V_{q,\theta}$, respectively. We would like to integrate by parts with respect to the variable $\xi$. To this end, observe that
$$\nabla_\xi\varphi=\frac{t}{2\tau(t-\tau)}(\xi-\zeta),$$
where
$$\zeta=\frac{\tau}{t}x+\frac{t-\tau}{t}y.$$ 
We are going to use the identity
$$e^{i\varphi}=L_{\xi}e^{i\varphi},$$
where
$$L_{\xi}=-i\frac{\nabla_{\xi}\varphi}{|\nabla_{\xi}\varphi|^2}\cdot \nabla_{\xi}=-i\frac{2\tau(t-\tau)}{t}\Lambda_{\xi-\zeta},$$
$$\Lambda_\xi=\frac{\xi}{|\xi|^2}\cdot\nabla_\xi.$$
It is easy to see by induction that given any integer $m\ge 0$, the operator $\left(\Lambda_\xi^*\right)^m$ is of the form
$$\left(\Lambda_\xi^*\right)^m=\sum_{0\le|\alpha|\le m}a_\alpha^{(m)}(\xi)\partial_\xi^\alpha,\eqno{(4.8)}$$
with functions $a_\alpha^{(m)}\in C^\infty({\bf R}^n\setminus 0)$ satisfying
$$\left|\partial_\xi^\beta a_\alpha^{(m)}(\xi)\right|\le C|\xi|^{-2m+|\alpha|-|\beta|},\quad\forall\xi\in {\bf R}^n\setminus 0.\eqno{(4.9)}$$
Given an integer $0\le m<(n-1)/2$ and a function $W\in C^m_0({\bf R}^n)$, we have
$$\int_{{\bf R}^n}e^{i\varphi}W(\xi)d\xi=\int_{{\bf R}^n}e^{i\varphi}\left(L_\xi^*\right)^mW(\xi)d\xi$$
$$=(-2i)^mt^{-m}\tau^m(t-\tau)^m\sum_{0\le|\alpha|\le m}\int_{{\bf R}^n}e^{i\varphi}a_\alpha^{(m)}(\xi-\zeta)\partial_\xi^\alpha W(\xi)d\xi.$$
Therefore, making a change of variables $\mu=t/\tau$ we can write
$$J_p(x,y,t,\gamma;W):=\int_0^{\gamma}\int_{{\bf R}^n}e^{i\varphi}(t-\tau)^{-n/2}\tau^{-n/2}\phi_p(\gamma/\tau)W(\xi)d\xi d\tau$$ 
$$=(-2i)^mt^{m-n+1}\sum_{0\le|\alpha|\le m}\int_{t/\gamma}^\infty\int_{{\bf R}^n}e^{i\varphi}\left(\frac{\mu}{\mu-1}\right)^{m-n/2}\mu^{n/2-2-m}\phi_p(\gamma\mu/t)$$ $$\times\, a_\alpha^{(m)}\left(\xi-y-\mu^{-1}(x-y)\right)\partial_\xi^\alpha W(\xi)d\xi d\mu.\eqno{(4.10)}$$
Setting
$$f(\mu)=\int_2^\mu \exp\left(i\frac{\lambda|x-\xi|^2}{4(\lambda-1)t}+i\frac{\lambda|y-\xi|^2}{4t}\right)\lambda^{-1/2}d\lambda,$$
$$g_\alpha^{(m)}(\mu)=\left(\frac{\mu}{\mu-1}\right)^{m-n/2}\mu^{(n-3)/2-m}\phi_p(\gamma\mu/t)a_\alpha^{(m)}\left(\xi-y-\mu^{-1}(x-y)\right),$$
we can write the integral in the RHS of (4.10) as follows
$$\int_{t/\gamma}^\infty\int_{{\bf R}^n}g_\alpha^{(m)}(\mu)\partial_\xi^\alpha W(\xi)d\xi df(\mu)=\int_{{\bf R}^n}f(\mu)g_\alpha^{(m)}(\mu)\left|_{t/\gamma}^\infty\right.\partial_\xi^\alpha W(\xi)d\xi$$ 
$$-\int_{t/\gamma}^\infty\int_{{\bf R}^n}f(\mu)\partial_\mu g_\alpha^{(m)}(\mu)\partial_\xi^\alpha W(\xi)d\xi d\mu.\eqno{(4.11)}$$
In view of (4.9), we have
$$\left|g_\alpha^{(m)}(\mu)\right|\le C\left(2^pt/\gamma\right)^{(n-3)/2-m}\left|\xi-y-\mu^{-1}(x-y)\right|^{-2m+|\alpha|},\eqno{(4.12)}$$
$$\left|\partial_\mu g_\alpha^{(m)}(\mu)\right|\le C\left(2^pt/\gamma\right)^{(n-3)/2-m-1}$$ $$\times\left(\left|\xi-y-\mu^{-1}(x-y)\right|^{-2m+|\alpha|}+\mu^{-1}|x-y|\left|\xi-y-\mu^{-1}(x-y)\right|^{-2m+|\alpha|-1}\right)$$ $$\le C\left(2^pt/\gamma\right)^{(n-3)/2-m-1}$$ $$\times\left(\left|\xi-y-\mu^{-1}(x-y)\right|^{-2m+|\alpha|}+|y-\xi|\left|\xi-y-\mu^{-1}(x-y)\right|^{-2m+|\alpha|-1}\right).\eqno{(4.13)}$$
We will also need the following

\begin{lemma} There exists a constant $C>0$ such that
$$|f(\mu)|\le Ct^{1/2}|y-\xi|^{-1},\quad\forall\mu\ge 2.\eqno{(4.14)}$$
\end{lemma}

By (4.10)-(4.14), we get
$$\left|J_p(x,y,t,\gamma;W)\right|\le Ct^{-n/2}\left(2^p/\gamma\right)^{(n-3)/2-m}\sum_{0\le|\alpha|\le m}\int_{{\bf R}^n}\left(\left|\xi-y\right|^{-2m+|\alpha|-1}\right.$$ $$\left.+|y-\xi|^{-1}\left|\xi-y-\gamma t^{-1}(x-y)\right|^{-2m+|\alpha|}+\left|\xi-y-\gamma t^{-1}(x-y)\right|^{-2m+|\alpha|-1}\right)\left|\partial_\xi^\alpha W(\xi)\right|d\xi$$
 $$+Ct^{-n/2-1}\left(2^p/\gamma\right)^{(n-3)/2-m-1}\sum_{0\le|\alpha|\le m}\int_{2^{p-1}t/\gamma}^{2^{p}t/\gamma}\int_{{\bf R}^n}\left(|y-\xi|^{-1}\left|\xi-y-\mu^{-1}(x-y)\right|^{-2m+|\alpha|}\right.$$ $$\left.+\left|\xi-y-\mu^{-1} (x-y)\right|^{-2m+|\alpha|-1}\right)\left|\partial_\xi^\alpha W(\xi)\right|d\xi d\mu$$ 
 $$\le Ct^{-n/2}\left(2^p/\gamma\right)^{(n-3)/2-m}\sum_{0\le|\alpha|\le m}\sup_{\xi\in{\bf R}^n}\langle\xi\rangle^{n+\epsilon-2m-1+|\alpha|}\left|\partial_\xi^\alpha W(\xi)\right|,\eqno{(4.15)}$$
for every $0<\epsilon\ll 1$. Applying (4.15) with $m=k_0$, $W=V_q-V_{q,\theta}$, and using (4.7), we obtain
$$\left|K_{p,q}^{(1)}\right|\le Ct^{-n/2}\theta^\nu\left(\gamma 2^{-p}\right)^{k_0-(n-3)/2}2^{-q(2k_0-n+3)}. \eqno{(4.16)}$$
We would like to apply (4.15) with $m=k_0+1$ as well. When $n$ is odd, however, this is impossible because (4.15) does not hold with 
$m=(n-1)/2$ (due to the fact that the function $|\xi|^{-n}$ is not integrable at $\xi=0$). Therefore, we need to make some modifications in this case. Consider first the case $n$ even. Then we can apply (4.15) with $m=k_0+1=(n-2)/2$, $W=V_{q,\theta}$, and use (4.5) and (4.6) to obtain
$$\left|K_{p,q}^{(2)}\right|\le Ct^{-n/2}\theta^{-1+\nu}\left(\gamma 2^{-p}\right)^{k_0-(n-5)/2}2^{-q(2k_0-n+5)}. \eqno{(4.17)}$$
By (4.16) and (4.17) we conclude
$$\left|K_{p,q}\right|\le Ct^{-n/2}\theta^\nu\left(\gamma 2^{-p}\right)^{k_0-(n-3)/2}2^{-q(2k_0-n+3)}\left(1+\theta^{-1}\gamma 2^{-p}2^{-2q}\right). \eqno{(4.18)}$$
Taking $\theta=\gamma 2^{-p}2^{-2q}$ we deduce from (4.18)
$$\left|K_{p,q}\right|\le Ct^{-n/2}\left(\gamma 2^{-p}\right)^{k-(n-3)/2}2^{-q(2k-n+3)}, \eqno{(4.19)}$$
which clearly implies (4.2) in this case. Let now $n$ be odd. Given an integer $0\le m<n/2$, by (4.9) and (4.10),
$$\left|J_p(x,y,t,\gamma;W)\right|\le Ct^{-n/2-1}\left(2^p/\gamma\right)^{(n-2)/2-m-1}$$ $$\times\sum_{0\le|\alpha|\le m}\int_{2^{p-1}t/\gamma}^{2^{p}t/\gamma}\int_{{\bf R}^n}\left|\xi-y-\mu^{-1} (x-y)\right|^{-2m+|\alpha|}\left|\partial_\xi^\alpha W(\xi)\right|d\xi d\mu$$ 
 $$\le Ct^{-n/2}\left(2^p/\gamma\right)^{(n-2)/2-m}\sum_{0\le|\alpha|\le m}\sup_{\xi\in{\bf R}^n}\langle\xi\rangle^{n+\epsilon-2m+|\alpha|}\left|\partial_\xi^\alpha W(\xi)\right|,\eqno{(4.20)}$$
for every $0<\epsilon\ll 1$. Applying (4.20) with $m=k_0+1=(n-1)/2$, $W=V_{q,\theta}$, and using (4.5) and (4.6), we get
$$\left|K_{p,q}^{(2)}\right|\le Ct^{-n/2}\theta^{-1+\nu}\left(\gamma 2^{-p}\right)^{1/2}2^{-q}. \eqno{(4.21)}$$
By (4.16) and (4.21) we conclude
$$\left|K_{p,q}\right|\le Ct^{-n/2}\theta^\nu\left(1+\theta^{-1}\left(\gamma 2^{-p}\right)^{1/2}2^{-q}\right). \eqno{(4.22)}$$
Taking $\theta=\left(\gamma 2^{-p}\right)^{1/2}2^{-q}$ we deduce from (4.22)
$$\left|K_{p,q}\right|\le Ct^{-n/2}\left(\gamma 2^{-p}\right)^{\nu/2}2^{-\nu q}, \eqno{(4.23)}$$
which implies (4.2) in this case. 
\eproof

{\it Proof of Lemma 4.2.} We need to show that the function
$$\widetilde f(\mu)=\int_2^\mu e^{i\varphi(\lambda)}\lambda^{-1/2}d\lambda,$$
where
$$\varphi(\lambda)=\frac{\lambda\sigma_1}{\lambda-1}+\lambda\sigma_2,$$
satisfies the bound
$$\left|\widetilde f(\mu)\right|\le C\sigma_2^{-1/2},\quad \forall\mu\ge 2,\,\sigma_1,\sigma_2>0,\eqno{(4.24)}$$
with a constant $C>0$ independent of $\mu,\,\sigma_1$ and $\sigma_2$. We have
$$\varphi':=\frac{d\varphi}{d\lambda}=\sigma_2-\sigma_1(\lambda-1)^{-2},$$
so $\varphi'$ vanishes at $\lambda_0=1+(\sigma_1/\sigma_2)^{1/2}$. Suppose first that $\lambda_0\ge 20/11$. We will consider several cases.

Case 1. $2\le\mu\le 9\lambda_0/10$. We have, for $2\le\lambda\le 9\lambda_0/10$,
$$\left|\varphi'(\lambda)\right|\ge \frac{\sigma_2|\lambda-\lambda_0|}{\lambda-1}\ge\frac{\sigma_2\lambda_0}{10(\lambda-1)},$$
$$\left|\frac{\varphi''(\lambda)}{\varphi'(\lambda)}\right|\le \frac{2}{\lambda-1}\left(1+\frac{\sigma_2}{\left|\varphi'(\lambda)\right|}\right)\le\frac{22}{\lambda-1}.$$
Integrating by parts, we get
$$\widetilde f(\mu)=\int_2^\mu (i\varphi')^{-1}\lambda^{-1/2}de^{i\varphi}=e^{i\varphi}(i\varphi')^{-1}\lambda^{-1/2}\left|_2^\mu-
\int_2^\mu e^{i\varphi}g(\lambda)d\lambda\right.,$$
where the function
$$g(\lambda)=\frac{d}{d\lambda}\left((i\varphi')^{-1}\lambda^{-1/2}\right)$$
satisfies the bound
$$|g(\lambda)|\le C\sigma_2^{-1}\lambda_0^{-1}\lambda^{-1/2}.$$
Hence,
$$|\widetilde f(\mu)|\le C\sigma_2^{-1}\lambda_0^{-1}\mu^{1/2}.\eqno{(4.25)}$$
On the other hand, we have
$$|\widetilde f(\mu)|\le 2\mu^{1/2}.\eqno{(4.26)}$$
By (4.25) and (4.26),
$$|\widetilde f(\mu)|\le C\sigma_2^{-1/2}\left(\mu/\lambda_0\right)^{1/2}\le C\sigma_2^{-1/2}.\eqno{(4.27)}$$

Case 2. $9\lambda_0/10\le\mu\le 11\lambda_0/10$. Making a change of variables $\lambda=\lambda_0(1+z)$ we can write
$$\widetilde f_1(\mu):=\int_{9\lambda_0/10}^\mu e^{i\varphi(\lambda)}\lambda^{-1/2}d\lambda=\lambda_0^{1/2}\int_{I}e^{i\kappa\phi(z)}(1+z)^{-1/2}dz,$$
where $\kappa=\sigma_2\lambda_0$, $I=[-10^{-1},\mu/\lambda_0-1]\subset [-10^{-1},10^{-1}]$ and 
$$\phi(z)=(1+z)\left(1+\frac{(\lambda_0-1)^2}{\lambda_0(1+z)-1}\right)=\lambda_0+\frac{\lambda_0}{\lambda_0-1}z^2+O(z^3),\quad |z|\le 10^{-1},$$
uniformly in $\lambda_0$. If $|a|\le 10^{-1}$ and $0<\theta_0\ll 1$, we can change the contour of integration to obtain
$$\left|\int_0^ae^{i\kappa\phi(z)}(1+z)^{-1/2}dz\right|$$ $$\le \left|\int_0^ae^{i\kappa\phi(e^{i\theta_0}y)}(1+e^{i\theta_0}y)^{-1/2}dy\right|
+\left|a\int_0^{\theta_0}e^{i\kappa\phi(e^{i\theta}a)}(1+e^{i\theta}a)^{-1/2}d\theta\right|$$
$$\le C'_1\int_0^a e^{-C_1\kappa y^2}dy+C'_2\int_0^{\theta_0}e^{-C_2\kappa\theta}d\theta=O(\kappa^{-1/2}).$$
Hence
$$|\widetilde f_1(\mu)|\le  C\sigma_2^{-1/2}.\eqno{(4.28)}$$

Case 3. $\mu> 11\lambda_0/10$. Changing the contour of integration we can write 
$$\widetilde f_2(\mu):=\int_{\mu_0}^\mu e^{i\varphi(\lambda)}\lambda^{-1/2}d\lambda=\int_0^{z_0}e^{i\varphi(\mu_0+ze^{i\theta_0})}
\left(\mu_0+ze^{i\theta_0}\right)^{-1/2}dz$$
$$+z_0\int_0^{\theta_0}e^{i\varphi(\mu_0+z_0e^{i\theta})}
\left(\mu_0+z_0e^{i\theta}\right)^{-1/2}d\theta,$$
where $0<\theta_0\ll 1$, $\mu_0=11\lambda_0/10$, $z_0=\mu-11\lambda_0/10$. It is easy to check that
$${\rm Im}\,\varphi(\mu_0+ze^{i\theta})\ge Cz\sigma_2\theta,\quad\forall z\ge 0,\, 0<\theta\le\theta_0,$$
with a constant $C>0$ independent of $z$, $\theta$ and $\sigma_2$. Therefore, we get
$$|\widetilde f_2(\mu)|\le \int_0^\infty e^{-C'z\sigma_2}z^{-1/2}dz+z_0^{1/2}\int_0^{\theta_0}e^{-Cz_0\sigma_2\theta}d\theta$$
 $$\le O\left(\sigma_2^{-1/2}\right)+z_0^{1/2}O\left((\sigma_2z_0)^{-1/2}\right)\le C\sigma_2^{-1/2}.\eqno{(4.29)}$$
 Clearly, when $\lambda_0\ge 20/11$ the bound (4.24) follows from (4.27)-(4.29). When $\lambda_0\le 20/11$ we have
 $${\rm Im}\,\varphi(2+ze^{i\theta})\ge Cz\sigma_2\theta,\quad\forall z\ge 0,\, 0<\theta\le\theta_0.$$
Therefore, proceeding as in Case 3 we conclude that (4.24) holds in this case, too.
\eproof

\section{Proof of Proposition 3.5}

We keep the same notations as in the previous section. Clearly, it suffices to prove the following

\begin{prop} Under the assumptions of Proposition 3.5, for every $0<\epsilon\ll 1$ there exists a constant $C_\epsilon>0$ so that we have the bound
$$\left|K_{p,q}(x,y,t,\gamma)\right|\le C_\epsilon \left\|V\right\|_{{\cal V}_\delta^{(n-1)/2}}2^{-\epsilon(p+q)}\gamma^{1-\epsilon} t^{-n/2},\quad t\ge 2\gamma.\eqno{(5.1)}$$
\end{prop}

{\it Proof.} Since $V\in {\cal V}^{(n-1)/2}_\delta({\bf R}^n)$, we have
$$\left|\partial_\xi^\alpha V_q(\xi)\right|\le C\langle\xi\rangle^{-\delta-|\alpha|},\quad 0\le|\alpha|\le (n-1)/2,\eqno{(5.2)}$$
and $\langle\xi\rangle\sim 2^q$ on supp$\,V_q(\xi)$, where the constant $C$ is of the form
$$C=C'\left\|V\right\|_{{\cal V}_\delta^{(n-1)/2}},\eqno{(5.3)}$$
with $C'>0$ independent of $V$. 
Choose a function $\eta\in C_0^\infty({\bf R}^n)$, $\eta(\xi)=1$ for $|\xi|\le 1$, $\eta(\xi)=0$ for $|\xi|\ge 2$ and let 
$0<\theta\le 1$ be a parameter to be fixed later on. Decompose the function $K_{p,q}$ as $\widetilde K_{p,q}^{(1)}+\widetilde K_{p,q}^{(2)}$, where $\widetilde K_{p,q}^{(1)}$ and $\widetilde K_{p,q}^{(2)}$ are defined by replacing $V_q$ in the definition of $K_{p,q}$ by $V_q-\widetilde V_{q,\theta}$ and $\widetilde V_{q,\theta}$, respectively, where
$$\widetilde V_{q,\theta}(\xi)=\eta\left(\theta^{-1}(\xi-\zeta)\right)V_q(\xi).$$
In view of (5.2) we have (for $0\le |\alpha|\le (n-1)/2$)
$$\left|\partial_\xi^\alpha\left(V_q(\xi)-\widetilde V_{q,\theta}(\xi)\right)\right|\le C2^{-q}\theta^{-\epsilon}\sum_{j=0}^{|\alpha|}|\xi-\zeta|^{-j+\epsilon}\langle\xi\rangle^{j-|\alpha|-2\epsilon},$$
which implies 
$$\sum_{0\le |\alpha|\le (n-1)/2}\int_{{\bf R}^n}\left(|y-\xi|^{-1}|\xi-\zeta|^{|\alpha|+1-n}+|\xi-\zeta|^{|\alpha|-n}\right)
\left|\partial_\xi^\alpha\left(V_q(\xi)-\widetilde V_{q,\theta}(\xi)\right)\right|d\xi$$
$$\le C2^{-q}\theta^{-\epsilon}\sum_{0\le |\alpha|\le (n-1)/2}\sum_{j=0}^{|\alpha|}$$ $$\int_{{\bf R}^n}\left(|y-\xi|^{-1}|\xi-\zeta|^{|\alpha|+1-n-j+\epsilon}+|\xi-\zeta|^{|\alpha|-n-j+\epsilon}\right)\langle\xi\rangle^{j-|\alpha|-2\epsilon}d\xi\le C2^{-q}\theta^{-\epsilon}.$$
Therefore, proceeding as in (4.15) (with $m=(n-1)/2$), we get
$$\left|\widetilde K_{p,q}^{(1)}\right|\le Ct^{-n/2}\theta^{-\epsilon}\gamma 2^{-p-q}.\eqno{(5.4)}$$
Similarly, we have
$$\left|\partial_\xi^\alpha\widetilde V_{q,\theta}(\xi)\right|\le C\theta\sum_{j=0}^{|\alpha|}|\xi-\zeta|^{-j-1}\langle\xi\rangle^{j-1-|\alpha|-\epsilon},$$
which implies 
$$\sum_{0\le |\alpha|\le (n-3)/2}\int_{{\bf R}^n}\left(|y-\xi|^{-1}|\xi-\zeta|^{|\alpha|+3-n}+|\xi-\zeta|^{|\alpha|+2-n}\right)
\left|\partial_\xi^\alpha\widetilde V_{q,\theta}(\xi)\right|d\xi$$
$$\le C\theta\sum_{0\le |\alpha|\le (n-3)/2}\sum_{j=0}^{|\alpha|}$$ $$\int_{{\bf R}^n}\left(|y-\xi|^{-1}|\xi-\zeta|^{|\alpha|+2-n-j}+|\xi-\zeta|^{|\alpha|+1-n-j}\right)\langle\xi\rangle^{j-1-|\alpha|-\epsilon}d\xi\le C\theta.$$
Therefore, proceeding as in (4.15) (with $m=(n-3)/2$), we get
$$\left|\widetilde K_{p,q}^{(2)}\right|\le Ct^{-n/2}\theta.\eqno{(5.5)}$$
By (5.4) and (5.5), we conclude
$$\left|K_{p,q}\right|\le Ct^{-n/2}\left(\theta^{-\epsilon}\gamma 2^{-p-q}+\theta\right)\le Ct^{-n/2}\left(\gamma 2^{-p-q}\right)^{1-\epsilon},\eqno{(5.6)}$$
if we choose $\theta=\gamma 2^{-p-q}$, where the constant $C$ is of the form (5.3). 
\eproof

F. Cardoso

Universidade Federal de Pernambuco, 

Departamento de Matem\'atica, 

CEP. 50540-740 Recife-Pe, Brazil

e-mail: fernando@dmat.ufpe.br\\

C. Cuevas

Universidade Federal de Pernambuco, 

Departamento de Matem\'atica, 

CEP. 50540-740 Recife-Pe, Brazil

e-mail: cch@dmat.ufpe.br\\

G. Vodev

Universit\'e de Nantes,

 D\'epartement de Math\'ematiques, UMR 6629 du CNRS,
 
 2, rue de la Houssini\`ere, BP 92208, 
 
 44332 Nantes Cedex 03, France

e-mail: georgi.vodev@math.univ-nantes.fr

\end{document}